\documentclass[11pt]{article}
\usepackage{amssymb}
\usepackage{mathtools}
\usepackage{mathrsfs}
\usepackage{amsmath}
\usepackage{amsmath}
\usepackage{amsfonts}
\usepackage{stmaryrd}
\usepackage{enumitem,kantlipsum}

\usepackage{amsthm}
\usepackage{thmtools}
\usepackage{tikz-cd}

\input amssym.def
\input amssym
\input xypic
\usepackage{comment}
\def\demo{\noindent{\bf Proof. }}
\def\QED{\hfill$\Box$}
\newcommand{\rar}{\rightarrow}

\newtheorem*{Theorem*}{Theorem}
\newtheorem{Theorem}{Theorem}[section]
\newtheorem{Lemma}[Theorem]{Lemma}

\newtheorem{Proposition}[Theorem]{Proposition}
\newtheorem{Example}[Theorem]{Example}

\newtheorem{rem}[Theorem]{Remark}
\newtheorem{conv}[Theorem]{Convention}

\begin{document}
\title{
{{ \bf Existence of birational small Cohen-Macaulay modules over biquadratic extensions in mixed characteristic}}}

\author{ Prashanth Sridhar
\\
\begin{tabular}{c}
\vspace{0.2in} \      \\
Department of Mathematics, University of Kansas \\
Lawrence, Kansas 66045 \\
e-mail: prashanth@ku.edu
\end{tabular}}

\date{ }

\maketitle \vspace{-0.2in}

\begin{abstract}
\noindent   Let $S$ be an unramified regular local ring of mixed characteristic two and $R$ the integral closure 
of $S$ in a biquadratic extension of its quotient field obtained by adjoining roots of sufficiently general square free elements $f,g\in S$. Let $S^2$ denote the subring of $S$ obtained by lifting to $S$ the image of the Frobenius map on $S/2S$. When at least one of $f,g\in S^2$, we characterize the Cohen-Macaulayness of $R$ and show that $R$ admits a birational small Cohen-Macaulay module. It is noted that $R$ is not automatically Cohen-Macaulay in case $f,g\in S^2$ or if $f,g\notin S^2$.
\end{abstract}

\section{Introduction}
Small Cohen-Macaulay (CM) modules or maximal Cohen-Macaulay modules (see \ref{conv1}) are very classical objects that have been studied extensively in the literature. However, their existence over local rings that are not Cohen-Macaulay is far from clear. There are examples of non existence over local rings that are not catenary, but one may view this as pathological. In fact, a domain that admits a small CM module has to be universally catenary, see \cite{10.1007/BFb0068925}. Hochster conjectured that every complete local domain admits a small CM module, but a positive answer to this question is known only in very few cases. Note that the conjecture reduces to the integral closure of a complete regular local ring in a finite normal extension of its fraction field. So it is natural to first look at the case of extensions of the fraction field of a regular local ring with a ``nice" Galois group.
\\
Let $S$ be an unramified regular local ring and $L$ its quotient field. Let $K$ be a finite extension of $L$ and $R$ the integral closure of $S$ in $K$. In \cite{RO}, Roberts showed that if $K/L$ is Abelian and $[K:L]$ is not divisible by the characteristic of the residue field of $S$, then $R$ is Cohen-Macaulay. In particular, this applies to the equi-characteristic zero case. In the mixed characteristic $p$ scenario, the conclusion fails as shown in \cite{KO} and \cite{KA}. The examples in these articles were obtained by considering extensions adjoining a $p$-th root of a single element that is not square free. However, in case we adjoin a $p$-th root of a single square free element, the integral closure is Cohen-Macaulay and this fact is relatively easy to show. In contrast, as we shall see, the integral closure need not be Cohen-Macaulay in a finite square free tower of $p$-th roots. In fact, it fails to be Cohen-Macaulay even if we adjoin $p$-th roots of two square free elements.
\par
Now assume additionally that $S$ has mixed characteristic $p>0$. As discussed above, radical extensions of $L$ obtained by adjoining $n$-th roots of elements of $S$ with the property that $p$ divides $n$ are prime examples of the failure of Roberts's theorem in mixed characteristic. Moreover, the importance of radical extensions stem from Kummer theory, which says that Abelian extensions are repeated radical extensions under the presence of ``suitable roots of unity". In this paper and \cite{sridhar2021existence}, we investigate such extensions to understand the obstructions one faces when $p$ divides $n$, with a focus on the question of existence of birational small CM modules over the integral closure $R$. Note that if $p\nmid n$ or if it were the case that $S$ contained the rational numbers, the integral closure in repeated radical extensions is Cohen-Macaulay, see \ref{charzeroradical}.    
\par In \cite{KA} it is shown that the integral closure of $S$ in an extension of its quotient field obtained by adjoining the $p$-th root of an arbitrary element of $S$ admits a birational small CM module. In \cite{KATZ2021350}, under certain circumstances $R$ is known to admit a birational small CM module in extensions obtained by adjoining the $p^n$-th root of a single element. In this article, we investigate the question of existence of small Cohen-Macaulay modules over the integral closure $R$ of $S$ in an extension $K/L$ of degree $p^2$ obtained by adjoining $p$-th roots of sufficiently general square free elements $f,g\in S$. Roughly speaking, one may think of this as the case where $Gal(K/L)=\mathbb{Z}_p\times \mathbb{Z}_p$.
\par 
To this end, we can reduce to the case $f,g\in S^p$ when $S$ is complete with perfect residue field, where $S^p$ is the subring of $S$ obtained by lifting the Frobenius map on $S/pS$ to $S$, see \cite{DKPS}. However, this does not mean $R$ is Cohen-Macaulay when $f,g\notin S^p$ as we show in \ref{example}. But we note that the example admits a small CM algebra (\ref{conv1}) in \ref{smallCMalgebra}. Finally, the square free condition on the elements is natural since such an extension corresponding to arbitrary $f,g\in S$ can be embedded in a sufficiently large finite tower obtained by adjoining square free $p$-th roots. 
\par
In the preliminary section, we  prove some results in mixed characteristic that are independent of the choice of prime integer $p$. In Sections $3$ and $4$, we work in mixed characteristic two. In Theorem \ref{P11}, when $f,g\in S^2$, we characterize when $R$ is Cohen-Macaulay and show that in case $R$ itself is not Cohen-Macaulay, $R$ admits a birational small Cohen-Macaulay module. For a similar discussion that is independent of mixed characteristic, see \cite{sridhar2021existence}. However, the results there are not as complete as in the mixed characteristic two case. 
\par More precisely, let $f,g \in S$ be relatively prime elements that are either units but not squares in $S$ or squarefree non-units. Let $\omega^2=f$ and $\mu^2=g$. For $n,k\geq 1$ integers, let $S^{p^k\wedge p^n}\subset S$ be the multiplicative subset of $S$ consisting of elements expressible in the form $x^{p^k}+y\cdot p^n$ for some $x,y\in S$. The main result of this paper is
 \begin{Theorem*}\noindent$\mathbf{4.10.}$
Let $S$ be an unramified regular local ring of mixed characteristic two and $f,g\in S^2$.
\begin{enumerate}
\item $R$ is Cohen-Macaulay if and only if one of the following happens
\begin{enumerate}
    \item At least one of $S[\omega],S[\mu]$ is not integrally closed.
    \item $S[\omega],S[\mu]$ are both integrally closed and $fg\notin S^{2\wedge 4}$.
    \item $S[\omega],S[\mu]$ are both integrally closed, $fg\in S^{2\wedge4}$ and $\mathscr{I}:=(2,f,g)\subset S$ is a two generated ideal or all of $R$.
\end{enumerate}
\item If $R$ is not Cohen-Macaulay, $R$ admits a birational small CM module. 
\end{enumerate}
\end{Theorem*}
Sections $3$ and $4$ are essentially one long proof of \ref{P11}. In section $3$ we prove the sufficiency of conditions $(a)$ and $(b)$ in \ref{T1} and \ref{T1'} respectively. Then in section $4$, we complete the proof. Surprisingly, the non Cohen-Macaulay cases occur only when the hypersurfaces $S[\omega]$ and $S[\mu]$ are both integrally closed.
\par For a ring $A$ and an $A$-module $M$, let $M^*$ denote the dual module $Hom_A(M,A)$. The strategy for proving \ref{P11} is to choose suitable ideals $J,I\subseteq A:=S[\omega,\mu]$ such that $R=J^*$ and $I^*$ is a $J^*$ module that is $S$-free.
We do not require in general the full strength of an unramified regular local ring $S$ for our results. In fact, many of the conclusions hold if we only assume that $S$ is a Noetherian integrally closed domain, $p\in S$ is prime and/or $S/pS$ is integrally closed.
\section{Preliminaries}
In this section we fix notation and record some observations that will be used subsequently.
Throughout this paper all rings considered are commutative and Noetherian.
\begin{conv}\normalfont\label{conv1}
\noindent \begin{itemize}
\item Let $(R,\mathfrak{m})$ be a local ring of dimension $d$. A nonzero $R$-module $M$ is a \textbf{small CM module} if it is finitely generated and every system of parameters of $R$ (equivalently some system of parameters) is a regular sequence on $M$. If $R$ is an arbitrary Noetherian ring, then an $R$-module $M$ is a \textbf{small CM module} if for all maximal ideals $\mathfrak{m}\subseteq R$, $M_{\mathfrak{m}}$ is a small CM module over $R_{\mathfrak{m}}$.
\item A Noetherian ring $R$ admits a \textbf{small CM algebra} $S$ if there is an injective, module finite map of rings $R\rightarrow S$ such that $S$ is Cohen-Macaulay.
\item For $R$ any commutative ring and $M$ an $R$-module, we will denote by $M^*_R$, the dual module $Hom_R(M,R)$. If $R$ is clear from the context, we will simply denote it by $M^*$. In particular, if $R$ is a domain with field of fractions $K$ and $M\subseteq K$, we use $M^*$ to also denote $(R:_K M)$ (see \cite{HS}[2.4.2] for example).
\item For a Noetherian ring $R$ of dimension at least one, we use the notation \[NNL_1(R):=\{P\in Spec(R)\: | \: height(P)=1,\: R_P\text{ is not a DVR}\}\]
\item Suppose $S$ is a ring and let $p\in \mathbb{Z}$ be a prime such that $p\in S$ is a non-unit. Let $F:S/pS\rightarrow S/pS$ be the Frobenius map. Let $S^p$ denote the subring of $S$ obtained by lifting the image of $F$ to $S$. Define $S^{p^k\wedge p^n}$ for $k,n\geq 1$ to be the multiplicative subset of $S$ of elements expressible in the form $x^{p^k}+y\cdot p^n$ for some $x,y\in S$. In particular, $S^{p\wedge p}=S^p$.
\item For a Noetherian local ring $R$ and $M$ a finite $R$-module, denote by $Syz^i_R(M)$ the $i$-th syzygy of $M$ in a minimal free resolution of $M$ over $R$.
\end{itemize}
\end{conv}
\begin{rem}\label{equicharzero}\normalfont
Suppose that $R\rightarrow T$ is a module finite extension of domains with $R$ integrally closed. Set $d:=[Frac(T):Frac(R)]$ and suppose $d\in R$ is a unit. Then if $T$ is Cohen-Macaulay, so is $R$. To see this, note that the trace map of their corresponding fraction fields gives an $R$-linear retraction $d^{-1}Tr:T\rightarrow R$. This ensures that for every ideal $I\subseteq R$, $IT\cap R=I$. Hence if $R$ is not Cohen-Macaulay, $T$ is not Cohen-Macaulay.

\end{rem}
\begin{rem}\label{genfact}\normalfont
 Recall the general fact: Let $R$ be a domain with field of fractions $L$ and let $K$ be a finite field extension of $L$. Then if the monic minimal polynomial $f(X)$ of $\gamma\in K$ over $L$ is such that $f(X)\in R[X]$, then $R[\gamma]\simeq R[X]/(f(X))$.
\end{rem}

\begin{Proposition}\label{charzeroradical}
Let $S$ be a regular local ring and $L$ its fraction field.
\\ Let $K:=L(\sqrt[n_1]{a_1},\dots,\sqrt[n_k]{a_k})$ with $a_i\in S$ and the $n_i$ positive integers that are units in $S$ for $1\leq i\leq k$. Then the integral closure of $S$ in $K$ is Cohen-Macaulay. In particular
\begin{enumerate}
\item If $S$ contains the rational numbers, the above conclusion holds for arbitrary integers $n_i$.
\item If $S$ has mixed characteristic $p>0$, the above conclusion holds for integers $n_i$ with the property that $p\not| \:n_i$ for all $1\leq i\leq k$.
\end{enumerate}
\end{Proposition}
\demo
First assume that the $a_i\in S$ are square free, mutually coprime and $n=n_i$ for all $1\leq i\leq k$. Then, from \cite{HK}[Prop 5.3], $R=S[\sqrt[n]{a_1},\dots,\sqrt[n]{a_k}]$ is integrally closed. Moreover, $R$ is Cohen-Macaulay from \ref{genfact}. So the conclusion holds in this case.
\par 
Suppose the $a_i$ and $n_i$ are arbitrary. Set $n:=\prod_{i=1}^k n_i$. Since $S$ is a UFD, there exist $b_1,\dots,b_m\in S$ square free and mutually coprime such that $K\hookrightarrow \mathscr{K}:=L(\sqrt[n]{b_1},\dots,\sqrt[n]{b_m})$. But the integral closure of $S$ in $\mathscr{K}$ is Cohen-Macaulay as observed above. By \ref{equicharzero}, the integral closure of $S$ in $K$ is Cohen-Macaulay.
\\ 
(1) and (2) now follow immediately.
\QED
\begin{conv}\normalfont\label{mainconv}
We will use the following notation for the remainder of the paper.
Let $S$ denote a Noetherian integrally closed domain of dimension $d$ and $L$ its field of fractions. Assume $Char(L)=0$. Let $p\in \mathbb{Z}$ be a prime such that $p\in S$ is a principal prime and $S/pS$ is integrally closed. An unramified regular local ring of mixed characteristic $p$ satisfies the above hypothesis, though not all results in this paper require this specific setting. The assumptions stated above will stand throughout the paper unless otherwise specified.
\par
An element $x\in S$ is said to be square free if for all height one primes $Q\subset S$ containing $x$, $QS_{Q}=(x)S_{Q}$. Say that a subset $W\subset S$ satisfies $\mathscr{A}_1$ if for all distinct $x,y\in W$, there exists no height one prime $Q\subset S$ such that $x,y\in Q$.
\par
Fix $f,g \in S$ such that they are not $p$-th powers in $S$, are square free and satisfy $\mathscr{A}_1$. Let $W,U$ be indeterminates over $S$. We have the monic irreducible polynomials $F(W):=W^p-f\in S[W]$ and $G(U):=U^p-g\in S[U]$. Let $K:=L(\omega,\mu)$ where $\omega$ and $\mu$ are roots of $F(W)$ and $G(U)$ respectively and assume that $G(U)$ is irreducible over $L(\omega)$, so that $[K:L]=p^2$. 
Denote by $R$ the integral closure of $S$ in $K$. That is, $R$ is the integral closure of $A:=S[\omega,\mu]$.
\end{conv}
\par
We make the following preliminary observations:

\begin{rem}\label{ACI}\normalfont
It follows from \ref{genfact} that $A\simeq S[W,U]/(F(W),G(U))$, $S[\omega]\simeq S[W]/(F(W))$ and $S[\mu]\simeq S[U]/(G(U))$.
\end{rem}

\begin{rem}\label{R1}\normalfont
We make use of the following observation later (\cite{V}[Theorem 2.4]): Let $S\subseteq C\subseteq D$ be an extension of Noetherian domains such that $S$ is integrally closed, $D$ is module finite over $S$ and $D$ is birational to $C$. Then if $C$ satisfies Serre's condition $R_1$, so does $D$. To see this, assume $C$ satisfies $R_1$. Let $P\subseteq D$ be any height one prime. Since going down holds for the extension $S\subseteq D$, $Q:=P\cap C$ is a height one prime in $C$. Since $C_Q\subseteq D_P$ is a birational extension and $C_Q$ is a DVR, we have $C_Q=D_P$. Thus $D$ satisfies $R_1$. 
\end{rem}

We include the following result from $\cite{KA}$ for convenience.
\begin{Proposition}\label{P2}
With notation as specified above, $S[\omega]$ is integrally closed if and only if $f\notin S^{p\wedge p^2}$. Further, if $S[\omega]$ is not integrally closed, write $f=h^p+a\cdot p^2$ for some $a,h\in S$, $h\neq 0$. Then
\begin{itemize}
\item[{\rm (a)}] $\overline{S[\omega]}=P^*_{S[\omega]}=S[\omega,\tau]$ where $\tau=p^{-1}\cdot(\omega^{p-1}+h\omega^{p-2}+\dots+h^{p-1})$ and $P:=(p,\omega-h)$ is the unique height one prime in $S[\omega]$ containing $p$.
\item[{\rm (b)}] If $p\geq 3$, there are exactly two height one primes in $\overline{S[\omega]}$ containing $p$, $Q_1$ and $Q_2$ satisfying ${Q_1}_{Q_1}=(\omega-h)_{Q_1}$ and ${Q_2}_{Q_2}=(p)_Q$. If $p=2$, $2\in S[\omega,\tau]=S[\tau]$ is square free.
\item[{\rm (c)}] If $p\geq 3$, $\tau$ satisfies $l(T):=T^2-cT-a\cdot(\omega-h)^{p-2}$ over $S[\omega]$ for some $c\in S[\omega]$. If $p=2$, $\tau$ satisfies $l(T):=T^2-hT-a$ over $S$.
\item[{\rm (d)}] $\overline{S[\omega]}$ is $S$-free with a basis given by the set $\{1,\omega,\dots,\omega^{p-2},\tau\}$.
\end{itemize}
\end{Proposition}

The next proposition characterizes when $A$ is integrally closed.
\begin{Proposition}\label{P1}
With established notation, the following hold:
\begin{enumerate}
    \item There exists a unique height one prime $P\subseteq A$ containing $p$.
    \item The ring $A$ is integrally closed if and only if $A_P$ is a DVR.
    \item The ring $A$ is integrally closed if and only if $f\notin S^p$, $g\notin S[\omega]_{(p)}^{p\wedge p^2}$ (or vice versa).
    \end{enumerate}
\end{Proposition}
\demo For (1), let $\phi:B:=S[W,U]\rightarrow A$ be the natural projection map. Height one primes in $A$ pull back to height three primes in $B$ containing $Ker(\phi)=(F(W),G(U))B$.
First assume that $f,g\in S^p$. Write $f=h_1^p+a\cdot p$ and $g=h_2^p+b\cdot p$ for some $h_1,h_2,a,b\in S$. It is then clear that the only height three prime in $B$ containing $Ker(\phi)$ and $p$ is $\tilde{P}\subset B$, given by $\tilde{P}:=(p,W-h_1,U-h_2)B$. Therefore $P:=(p,\omega-h_1,\mu-h_2)A$ is the unique height one prime in $A$ containing $p$ in this case.
\par 
Now let $f=h_1^p+a\cdot p$ and $g\notin S^p$. From \ref{P2}, $S[\mu]$ is integrally closed. Since $S/pS$ is integrally closed and $g\notin S^p$, $p\in S[\mu]$ is a principal prime. Since $A\simeq S[\mu][W]/(F(W))$, it is now clear that $P:=(p,\omega-h_1)A$ is the unique height one prime in $A$ containing $p$.
\par
Now assume that $f,g\notin S^p$. As noted above, $p$ is a principal prime in the integrally closed rings $S[\omega]$ and $S[\mu]$. We need to show that there exists a unique height three prime ideal of $B$ minimal over $(p,F(W),G(U))B$ or equivalently a unique height one prime minimal in $C[U]$ over $(G(U))C$ where $C:=S[\omega]/(pS[\omega])$ is a domain. Let $Q$ be the fraction field of $C$. If $G(U)$ is irreducible over $Q$, then from \ref{genfact} we get $C[\mu]\simeq C[U]/(G(U))$ so that $p$ is a principal prime in $A$. If $G(U)$ is reducible over $Q$, then $G(U)= (U-r)^p$ in $Q[U]$ for some $r\in Q$, so that $(U-r)Q[U]$ is the unique minimal prime over $G(U)Q[U]$. Since every prime $T\subseteq C[U]$ minimal over $G(U)C[U]$ intersects trivially with $C$, there is a unique height one prime $P\subseteq A$ containing $p$.
\par
For (2), observe that $A$ satisfies $S_2$ since it is $S$-free. To show the reverse implication, we see from part (1) that it suffices to show that $A[1/p]$ is integrally closed. Applying \cite{HK}[Proposition 5.3] to the ring $S[1/p]$, we see that $A[1/p]$ is indeed integrally closed. This completes the proof of the backward direction of (2). The forward direction is obvious.
\par 
For the backward direction of (3), assume that $f\notin S^p$ and $g\notin S[\omega]_{(p)}^{p\wedge p^2}$. As noted in the proof of part (1), $p\in S[\omega]$ is a prime. Suppose $g\notin S[\omega]_{(p)}^p$. Then $PA_P=pA_P$ and (2) implies that $A$ is integrally closed. Next, suppose that $g-h^p\in pS[\omega]_{(p)}$ for some $h\in S[\omega]_{(p)}$. Since $A_P\simeq S[\omega][U]_{(p,U-h)}/(G(U))$, $A_P$ is a DVR if and only if $G(U)\notin (p,U-h)^2S[\omega][U]_{(p,U-h)}$. Since \begin{equation}\label{cyclotomic}
U^{p-1}+\dots+h^{p-1}\in (p,U-h)S[\omega][U]
\end{equation} we see from our hypothesis that $G(U)\notin (p,U-h)^2S[\omega][U]_{(p,U-h)}$. Thus $A_P$ is a DVR and the conclusion follows from part (2).
\par 
For the forward direction of (3), we prove the contrapositive. As a first case, suppose $f,g\in S^p$ and $f=h_1^p+ap$, $g=h_2^p+bp$ with $h_1,h_2,a,b\in S$. For any $1\leq i\leq p-1$, notice that the element $\eta_i:=p^{-1}(\omega-h_1)^i(\mu-h_2)^{p-i}\in K$ satisfies $\eta_i^p\in A$ since $(\omega-h_1)^p,(\mu-h_2)^p\in pA$. But $\eta_i\notin A$ since $A$ is $S$-free with basis $\{\omega^i\mu^j\:|\: 0\leq i,j\leq p-1\}$. Thus $A$ is not integrally closed.
\par Now suppose that $f\notin S^p$ and $g\in S[\omega]_{(p)}^{p\wedge p^2}$. Let $g-h^p\in pS[\omega]_{(p)}$ for some $h\in S[\omega]_{(p)}$. We have $A_P\simeq S[\omega][U]_{(p,U-h)}/(G(U))$ since $p\in S[\omega]$ is prime. From our assumption and (\ref{cyclotomic}), it now follows that $A_P$ is not a DVR and hence $A$ is not integrally closed. This finishes the proof of the forward implication of (3).
\QED
\par
We first note a natural extension of \cite{KA}[3.2]:
\begin{Proposition}\label{P5}
With established notation, $R$ is $S$-free if $f\notin S^p$ and $g\in S[\omega]^p$. In particular, if $S$ is Cohen-Macaulay, then $R$ is Cohen-Macaulay.
\end{Proposition}
\demo
Since $f\notin S^p$, $S[\omega]$ is integrally closed by \ref{P2} and $p\in S[\omega]$ is a principal prime. Moreover,  \ref{P1}(3) allows us to assume that $g\in S[\omega]_{(p)}^{p\wedge p^2}$. Write $g=h^p+pb$, with $h,b\in S[\omega]$. Note that $g\in S[\omega]_{(p)}^{p\wedge p^2}$ implies that $b\in pS[\omega]_{(p)}\cap S=pS[\omega]$. That is $g\in S[\omega]^{p\wedge p^2}$. In this case, the proof of \cite{KA}[Lemma 3.2] goes through, so that $R$ is $S[\omega]$-free and hence $S$-free. Thus the proof is complete.
\QED
\\
\\
We need the following proposition for \ref{example}. The form given here is a bit more general than we actually need.
\begin{Proposition}\label{condlemma}
Let $T$ be any Gorenstein local domain such that its integral closure $T'$ is a finite $T$-module. Let $J$ denote the conductor ideal of $T$. Then $T'$ is Cohen-Macaulay if and only if $T/J$ is Cohen-Macaulay.
\end{Proposition}
\demo 
Let $E$ denote the field of fractions of $T$. Since $End(J):=(J:_E J)$ is a ring, we have $T'\subseteq J^*=End(J)\subseteq T'$, so that $J^*=T'$. For any $0\neq x\in J$, set $J':=(x:_T J)$. Note that $J$ is height one unmixed (see for example \cite{KA}[Proposition 2.1(2)]). Also, $J$ is principal if and only if $T$ is integrally closed, so we may assume $J$ is not principal. Since $T$ is Gorenstein, $T/J$ is Cohen-Macaulay if and only if $T/J'$ is Cohen-Macaulay (see \cite{10.2307/1971402}[Proposition 2.5] for example). From the depth lemma, $T/J'$ is a Cohen-Macaulay ring if and only if $J'\simeq T'$ is a Cohen-Macaulay $T$-module. This completes the proof since $T'$ is a Cohen-Macaulay ring if and only if it is a Cohen-Macaulay $T$-module.
\QED
\\
\\
When $f,g\notin S^p$, $R$ is not necessarily Cohen-Macaulay, as shown in the following example.
\begin{Example}\normalfont\label{example}
Set $S:=\mathbb{Z}[X,Y,V]_{(2,X,Y,V)}$. Let $f=XV^2+4$ and $g=XY^2+4$. Then $f,g$ are square free, form a regular sequence in $S$ and do not lie in $S^2$. Note that $2\in S[\omega]$ is a prime. Set $C:=S[\omega]/(2)\simeq (S/2S)[\gamma]$ where $\gamma=\sqrt{xv^2}$ and $x,y,v$ denote the respective images in $S/2S$. Since
$g=XY^2+4=(V^{-1}Y\omega)^2+4(1-V^{-2}Y^2)$, from \ref{P1}(3), $A$ is not integrally closed. Moreover from \cite{KA}[Lemma 3.2], $(P^*)_P=R_P$. From \ref{P1}(2), for all height one primes $Q\subseteq A$, $Q\neq P$, $(P^*)_Q=R_Q=A_Q$. Since $P^*$ and $R$ are birational $S_2$ $A$-modules, $P^*=R$. From \ref{P1}(2), the conductor of $A$ is contained in $P$ and hence is equal to $P$. We now show that $A/P$ is not Cohen-Macaulay, so that by \ref{condlemma}, $R$ is not Cohen-Macaulay.
\par Let $Q\subseteq C[U]$ denote the unique height one prime minimal over the image of $G(U)$. Set $\epsilon:=v^{-1}y\gamma$. Then $Q=(U-\epsilon)Frac(C)[U]\cap C[U]$ is the kernel of the natural surjection $C[U]\rightarrow C[\epsilon]$. Thus $A/P\simeq D:=C[\epsilon]$. 
\par Since $D$ is module finite over the regular local ring $S/2S$, it is Cohen-Macaulay if and only if it is $S/2S$-free. Certainly $D$ is generated over $S/2S$ by the set $\{ 1,\gamma,\epsilon,\gamma\epsilon\}$. Since $\epsilon \cdot \gamma=xyv\in S/2S$, we can trim this set to $G:=\{ 1,\gamma,\epsilon \}$. But $p.d._{S/2S}(D)=1$, since $D$ admits the minimal free resolution
 \begin{equation*}
\xymatrix@C+1pc{
 0\ar[r] & S/2S \ar[r]^{\psi^T}
  &
  (S/2S)^3  \ar[r]^{\phi}
  &
  D \ar[r] & 0 
}
\end{equation*}
where $\phi$ is the natural projection corresponding to the ordered set $G$ and $\psi=[ 0 \; y \; -v]$. Thus $A/P$ is not Cohen-Macaulay and hence $R$ is not Cohen-Macaulay.
\QED
\end{Example}
We will see in \ref{smallCMalgebra} that if $R$ is as in \ref{example}, then it admits a small CM algebra.
In general if $f,g\notin S^p$, we have not been able to construct a birational small CM module. However when $S$ is complete with perfect residue field, we can ``reduce" to the case $f,g\in S^p$ if we relax the birationality constraint, see \cite{DKPS}. Motivated by this, we focus on the case $f,g\in S^p$ for the remainder of the paper. 

\section{Sufficient conditions for $R$ to be Cohen-Macaulay}
In this section, we begin our proof of Theorem \ref{P11}. We shall demonstrate that conditions (a) and (b) are sufficient to give $R$ the Cohen-Macaulay property. In fact, we shall show that the sufficiency holds without the full hypothesis on $S$. To do this, we maintain the notation established in \ref{mainconv} and make the additional assumptions that $p=2$ and $f,g\in S^2$.
\\ 
Write $f=h_1^2+2\cdot a$ and $g=h_2^2+2\cdot b$ with $h_1,h_2,a,b\in S$. Note that under these assumptions, if $f\in 2S$, $f=2\cdot a$ for some $a\notin 2S$. This is because $f\in S$ is square free. So in this case $S[\omega]$ is necessarily integrally closed by \ref{P2}.
\begin{Proposition} \label{T1} $R$ is $S$-free if at least one of the rings $S[\omega],S[\mu]$ is not integrally closed.
\end{Proposition}
\demo 
 First assume that both $S[\omega]$ and $S[\mu]$ are not integrally closed. We have that $S[\tau_1]$ is integrally closed from $\ref{P2}$, for $\tau_1:=2^{-1}(\omega+h_1)$. Further $\tau_1$ satisfies $l_1(T):=T^2-h_1T-a'\in S[T]$ where $a=2a'$ for some $a'\in S$ and $T$ is an indeterminate over $S$. Further, $2\in S[\tau_1]$ is square free since $l_1(T)$ and $l_1'(T)$ are relatively prime over the quotient field of $S/2S$. Writing $b=2b'$ for some $b'\in S$, we also have that $l_2(T):=T^2-h_2T-b'$ and $l_2'(T)$ are relatively prime over the quotient field of $S[\tau_1]/Q$ for all height one primes $Q\subseteq S[\tau_1]$ containing $2$. Therefore $2\in E=S[\tau_1,\tau_2]$ is square free as well. Applying \cite{HK}[Proposition 5.3] to the ring $S[1/2]$, we see that $R[1/2]=A[1/2]\subseteq E[1/2]\subseteq R[1/2]$. Therefore $NNL_1(E)\subseteq V(2)$. Since $2\in E$ is square free, $E$ is regular in codimension one. Clearly $E$ is generated over $S$ by $\{1,\tau_1,\tau_2,\tau_1\tau_2 \}$ and hence $E$ is $S$-free of rank four. Thus $E$ satisfies Serre's criterion $S_2$ and is integrally closed, that is $E=R$.
 \par
 Next, without loss of generality assume $S[\mu]$ is integrally closed and $S[\omega]$ is not. From $\ref{P2}$, we have $S[\tau_1]$ is integrally closed for $\tau_1:=2^{-1}(\omega+h_1)$ and that $2\in S[\tau_1]$ is square free. Since $E:=S[\tau_1,\mu]\simeq S[\tau_1][U]/(G(U))$, height one primes in $E$ containing $2$ are of the form $(Q,\mu-h_2)E$ where $Q\subseteq S[\tau_1]$ is a height one prime containing $2$. By \ref{P1}(2) and \ref{R1}, $NNL_1(E)\subseteq V(2)$. But for any height one prime $\mathscr{P}:=(Q,\mu-h_2)\subseteq E$ containing $2$, $\mathscr{P}_{\mathscr{P}}=(\mu-h_2)_{\mathscr{P}}$. This is because going down holds for the extension $S\subseteq S[\mu]\subseteq E$, so $\mathscr{P}$ contracts back to the height one prime $P:=(2,\mu-h_2)\subseteq S[\mu]$. Since $(\mu-h_2)(\mu+h_2)=2b$ and $b\notin P$ by \ref{P2}, we have $PS[\mu]_P=(\mu-h_2)_P$. Thus, $\mathscr{P}_{\mathscr{P}}=(\mu-h_2)_{\mathscr{P}}$ and $E$ is regular in codimension one. Clearly $E$ is generated over $S$ by $\{ 1,\mu,\tau_1,\mu\tau_1\}$. Thus $E$ is $S$-free of rank four and hence satisfies Serre's criterion $S_2$. So $E=R$ and this completes the proof.
\QED
\\
\\
From \ref{equicharzero}, we see that a non Cohen-Macaulay normal domain containing the rationals does not admit a small CM algebra. In equal characteristic $p>0$, examples of non existence of small CM algebras are known, see \cite{BHATT20141}. As an immediate consequence of \ref{T1}, we record an example of the failure of this non-existence of small CM algebras in mixed characteristic.
\begin{Example}\label{smallCMalgebra}\normalfont
We assume notation as in \ref{example}, so that $R$ is a non Cohen-Macaulay normal domain of mixed characteristic $2$. Set $K':=L(\sqrt{X})$ and $T:=S[\sqrt{X}]$. Note that $T$ is an unramified regular local ring of mixed characteristic $2$ and $f,g\in T^{2\wedge 4}$. We claim that $f,g\in T$ are square free. To show this, we can assume that $2\in S$ is a unit since $f,g\notin 2S$. Then, by \cite{HK}[Proposition 5.2] $f,g\in T$ are square free. Clearly, $f,g\in T$ satisfy $\mathscr{A}_1$. Therefore by \ref{T1}, the integral closure of $T$ in $\mathscr{K}:=K'(\omega,\mu)$, say $\mathscr{R}$, is Cohen-Macaulay. Therefore, $\mathscr{R}$ is a small CM algebra over $R$.
\\
\begin{figure}[!h]
\centering
\begin{tikzcd}[row sep=60pt,column sep=huge]
  S \rar[hook] \dar[swap,hook]
  & T \dar[ hook] \rar[hook] & K' \dar[swap,hook,dd]\\
  R \rar[swap ,hook] \dar[swap ,hook]
  & \mathscr{R} \arrow[dr,hook,swap,sloped]\\
  K \rar[ rr, hook] && \mathscr{K}
\end{tikzcd}
\end{figure}
\end{Example}
\begin{Proposition}\label{T1'}
With established notation, $R$ is $S$-free if $S[\omega]$ and $S[\mu]$ are integrally closed and $fg\notin S^{2\wedge 4}$.
Further, in this case $P^*_A=R$ so that $P$ is the conductor of $R$ to $A$, where $P$ is the unique height one prime in $A$ containing $2$.
\end{Proposition}
\demo Since $S[\omega]$ and $S[\mu]$ are integrally closed, we have from $\ref{P2}$ that $f,g\notin S^{2\wedge 4}$. Write $f=h_1^2+2\cdot a$ and $g=h_2^2+2\cdot b$ with $a,b\notin 2S$. The condition $fg\notin S^{2\wedge 4}$ is equivalent to the condition $(ah_2^2+bh_1^2)\notin 2S$. This follows since
\begin{equation}\label{fgrel}
\begin{aligned}
fg&=(h_1^2+2a)(h_2^2+2b)
\\
&=(h_1h_2)^2+(ah_2^2+bh_1^2)\cdot 2+4ab
\end{aligned}
\end{equation} 
Note that the above is equivalent to requiring that $S[\omega\mu]$ be integrally closed. This is because, since $f,g\in S$ satisfy $\mathscr{A}_1$, $fg\in S$ is square free. Following this, $S[\omega\mu]$ is integrally closed if and only if $fg\notin S^{2\wedge 4}$ by \ref{P2}.

Let $\tau=2^{-1}(\mu-h_2)(\omega-h_1)\in K$. We see that $\tau$ satisfies \[l(T):= T^2-k_1k_2\in A[T]\] where $k_1:=2^{-1}(\omega-h_1)^2=h_1^2 + a -\omega h_1$, $k_2:=2^{-1}(\mu-h_2)^2=h_2^2 + b -\mu h_2$ and $T$ is an indeterminate over $A$. Note that
\begin{equation}\label{T1'42}
2=k_1^{-1}(\omega-h_1)^2=k_2^{-1}(\mu-h_2)^2
\end{equation} 
and $k_1,k_2\notin P$. We claim that $C:=S[\omega,\mu,\tau]$ is integrally closed under the given hypothesis. The unique height one prime $P\subseteq A$ containing $2$ is $P:=(2,\omega-h_1,\mu-h_2)$. Now \[l(T)\equiv T^2-ab \in (A/P)[T]\simeq (S/2S)[T]\] There exists a unique height one prime containing $2$ in $C$ and since $S/2S$ is integrally closed, the only possible forms for this unique height one prime are $Q_1:=PC$ or $Q_2:=(2, \omega-h_1,\mu-h_2,\tau-m)C$ for some $m\in S$ satisfying $m^2-ab\in 2S$. In the first case, ${Q_1}_{Q_1}$ is principal due to ($\ref{T1'42}$) and \begin{equation}\label{T1'41}
(\mu-h_2)\tau=k_2(\omega-h_1)
\end{equation}
Now let $Q_2\subseteq C$ be the unique height one prime in question. We will show that ${Q_2}_{Q_2}=(\tau-m)_{Q_2}$. First from (\ref{T1'42}) and (\ref{T1'41}), we have ${Q_2}_{Q_2}=(\tau-m,\mu-h_2)_{Q_2}$. By definition, we have in $C_{Q_2}$
\begin{equation}\label{T1'43}
\begin{aligned}
(\tau-m)^2&\equiv (k_1k_2+ab)\; \text{mod} (2) \\
&\equiv (ah_2(\mu-h_2)+bh_1(\omega-h_1)+h_1h_2(\omega-h_1)(\mu-h_2))\; \text{mod} (2)\\
&\equiv(\mu-h_2)(ah_2+bh_1k_2^{-1}\tau+h_1h_2(\omega-h_1))\; \text{mod} (2)
\end{aligned}
\end{equation}
where the last equivalence follows from (\ref{T1'41}). We claim that $ah_2+bh_1k_2^{-1}\tau$ is a unit in $C_{Q_2}$. If the claim holds, it follows from (\ref{T1'42}) that ${Q_2}_{Q_2}=(\tau-m)_{Q_2}$. To show the claim, assume on the contrary that $ah_2+bh_1k_2^{-1}\tau\in {Q_2}_{Q_2}$. Then \[k_2(a^2h_2^2+b^2h_1^2k_2^{-2}\tau^2)=k_2a^2h_2^2+b^2h_1^2k_1\in {Q_2}_{Q_2}\]
By definition of $k_1,k_2$ we get $(ah_2^2+bh_1^2)ab\in {Q_2}_{Q_2}$ and hence $(ah_2^2+bh_1^2)\in {Q_2}_{Q_2}\cap S=2S$. This contradicts our hypothesis. Thus the claim is true and ${Q_2}_{Q_2}$ is principal. From $\ref{P1}$(2) and \ref{R1}, $C$ is regular in codimension one. Let $D$ denote the $S$-module generated by $G:=\{1,\omega,\mu,\tau\}$. Note that $D$ is in fact a ring and is $S$-free of rank four. Then $D\subseteq C\subseteq D$. Thus $C=D$ satisfies $S_2$ and hence $C=R$ is $S$-free. 
\par We now show that $P$ is the conductor of $R$ to $A$. Since $A$ is not integrally closed, $A_P$ is not a DVR by \ref{P1}(2). Therefore the conductor is contained in $P$. On the other hand, since $P\cdot \tau\subseteq A$ and $R=A+S\cdot \tau$, $P$ conducts $R$ into $A$. Thus $P$ is the conductor of $R$ to $A$ and the proof is complete.
\QED
\section{Existence of birational small CM modules}

In this section we identify what it means for $R$ to be Cohen-Macaulay when $S$ is an unramified regular local ring of mixed characteristic $2$ and $f,g\in S^2$. When $R$ is not Cohen-Macaulay, we show the existence of a birational small CM module.
\par Towards this, from \ref{T1} and \ref{T1'}, if we seek a non $S$-free $R$, we must have that $S[\omega]$ and $S[\mu]$ are integrally closed such that $S[\omega\mu]\cong S[X]/(X^2-fg)$ is not integrally closed. This scenario is very much possible, see \ref{ex}. In this situation, we start by identifying an ideal $J\subseteq A$ such that $J^*=R$.
\begin{conv}\normalfont
For this section, we assume notation as specified at the beginning of section 3. In case we are in the situation $f,g\notin S^{2\wedge 4}$, assume that $f,g\notin 2S$. This is justified, since if exactly one of $f,g\in 2S$, then by $\ref{T1'}$, $R$ is $S$-free. The case $f,g\in 2S$ is not possible since $f,g$ satisfy $\mathscr{A}_1$.
\end{conv}
\begin{Proposition}\label{P10}
With established notation, let $S[\omega],S[\mu]$ be integrally closed and $fg\in S^{2\wedge 4}$, so that $S[\omega\mu]$ is not integrally closed. Then for $J:=(2,\omega\mu-h_1h_2)A$, we have $J^*_A=R$.
\end{Proposition}
\demo Since $J^*$ and $R$ are birational $S_2$ $A$-modules, it suffices to show $J^*_Q=R_Q$ for all height one primes $Q\subseteq A$. From \ref{P1}(2), $J^*_Q=R_Q=A_Q$ for all height one primes $Q\neq P$. So we may assume $(S,2S)$ and $(A,P)$ are one dimensional local rings.
\par Note that $A=S[\omega,\omega\mu]$ and that $\{f,fg\}$ satisfy $\mathscr{A}_1$ since they are both units. Since $S[\omega]$ is integrally closed and $S[\omega\mu]$ is not, the description of $R$ from the proof of \ref{T1} applies. We have that $R$ is generated over $S$ by the set $\{1,\mu,\tau,\mu\tau\}$ where $\tau=2^{-1}(\omega\mu+h_1h_2)$. This immediately implies $J$ conducts $R$ into $A$.
\par Let $\phi:B:=S[W,T]\rightarrow A$ be the projection map defined by $W\mapsto \omega$ and $T\mapsto \omega\mu$, where $W,T$ are indeterminates over $S$. Note that $Ker(\phi):=(W^2-f,T^2-fg)$. Suppose $l\in A$ conducts $R$ to $A$. Since $A_P$ is not regular, $l\in P=(2,\omega-h_1,\omega\mu-h_1h_2)$. Write $l=x\cdot 2+y\cdot (\omega-h_1)+z\cdot (\omega\mu-h_1h_2)$ for some $x,y,z\in A$. Viewing $l\cdot \tau\in A$ in $B$ and denoting lifts by $\sim$, we get
\begin{equation}
    \tilde{y}\cdot(W-h_1)(T-h_1h_2)\in (2, (W-h_1)^2,(T-h_1h_2)^2)
\end{equation}
By a standard regular sequence argument, $\tilde{y}\in (2,W-h_1,T-h_1h_2)$ and so $y\in P$. Since $(\omega-h_1)^2\in 2A$, we have $l\in J$. Thus $J$ is the conductor of $R$ to $A$ and the proof is complete.
\QED

\begin{Lemma}\label{Ifreeres}
With established notation, set \[I:=(2,\omega\mu-h_1h_2,h_2\omega-h_1\mu)A=(2,\omega\mu-h_1h_2,(\omega+h_1)(\mu+h_2))A\] Then $p.d._S(I)\leq 1$. More precisely, $I\simeq S^2\oplus_S C$ for some $S$-module $C$ that admits the free resolution
 \begin{equation*}
\xymatrix@C+1pc{
 0\ar[r] & S \ar[r]^{\psi^T}
  &
  S^3  \ar[r]^{\phi}
  &
  C \ar[r] & 0 
}
\end{equation*}
where $\phi$ is given by $\phi(e_1)=2\omega$, $\phi(e_2)=2\mu$ and $\phi(e_3)=h_2\omega-h_1\mu$ and $\psi=[ -h_2 \; h_1 \; 2]$.
\end{Lemma}
\demo 
We claim that $I$ is generated over $S$ by the set $G:=\{2,2\omega,2\mu,\omega\mu-h_1h_2,h_2\omega-h_1\mu\}$. To see this, note that $2\omega\mu\in  (\omega\mu-h_1h_2)\cdot S+2\cdot S$. Next $\omega\cdot (\omega\mu-h_1h_2)=a\cdot 2\mu-h_1(h_2\omega-h_1\mu)$. A symmetric argument takes care of $\mu\cdot (\omega\mu-h_1h_2)$. We also have $\omega\mu\cdot (\omega\mu-h_1h_2)=-h_1h_2(\omega\mu-h_1h_2)+4\cdot e$ for some $e\in S$. Finally, since $(\omega,\mu)\subseteq ((2,\omega\mu-h_1h_2):_A h_2\omega-h_1\mu)$, the claim holds.
\par Now, let $\phi':S^5\rightarrow I$ be the projection map defined by the ordered generating set $G$. If $[x_1\:x_2\:x_3\:x_4\:x_5]^T\in Ker(\phi)$, then since $A$ is $S$-free with a basis given by $\{1,\omega,\mu,\omega\mu\}$, we have that $x_1=x_4=0$. Therefore $I\simeq S^2\bigoplus_S C$, where $C$ is the $S$-module generated by $\{ 2\omega,2\mu,h_2\omega-h_1\mu\}$.  Now $C$ admits the above resolution for if $E:=[s_1\:s_2\:s_3]^T\in S^3$, then $E\in Ker(\phi)$ if and only if $2s_1+h_2s_3=0$ and $2s_2-h_1s_3=0$. Thus, if $E\in Ker(\phi)$, then there exists $k\in S$ such that $s_1=-h_2k$, $s_2=h_1k$ and $s_3=2k$, so that $E\in Im(\psi^T)$. \QED

\begin{rem}\label{Gor}\normalfont
We will use the following well known fact in the proposition below: Let $S\subseteq R$ be a finite extension of Noetherian local rings such that $S$ is Gorenstein and $R$ is Cohen-Macaulay. Then for any finite $R$-module $M$, $Hom_R(M,\omega_R)\simeq Hom_S(M,S)$ as $R$-modules (and $S$-modules), where $\omega_R$ is the canonical module of $R$. Indeed, we have $\omega_R\simeq Hom_S(R,S)$, so that by Hom-tensor adjointness, we have what we want.
\\
In our setting, let $(S,\mathfrak{m})$ be regular local, so that $(A,\mathfrak{n})$ is local where $\mathfrak{n}=(\mathfrak{m},\omega-h_1,\mu-h_2)$. Then for the extension $S\subseteq A$, since $A$ is Gorenstein, we have $Hom_A(M,A)\simeq Hom_S(M,S)$ as $A$-modules and $S$-modules.
\end{rem}

\begin{Proposition}\label{R2} 
Let $(S,\mathfrak{m})$ be an unramified regular local ring of mixed characteristic two. Let $S[\omega],S[\mu]$ be integrally closed rings and $fg\in S^{2\wedge 4}$. Then
\begin{enumerate}
\item If $f,g\in \mathfrak{m}$, $R\simeq S^2\bigoplus_S Syz^2_S(S/Q)$ where $Q:=(2,h_1,h_2)\subset S$.
\item $p.d._S(R)\leq 1$.
\item $R$ is Cohen-Macaulay if and only if $Q$ is a two generated ideal or all of $R$.
\end{enumerate}
\end{Proposition}
\demo We have from \ref{P10} that $J^*_A=R$, for $J:=(2,\omega\mu-h_1h_2)A$.  Let $I\subseteq A$ be as in \ref{Ifreeres} and $P:=(2,\omega-h_1,\mu-h_2)$ be the unique height one prime containing $2$ in $A$. Now $IA_P=JA_P$ since $\omega\notin P$ and $\omega\cdot (h_2\omega-h_1\mu)\in J$. Clearly $rad(I)=P$. Therefore, by \cite{KA}[Prop 2.1], $I^*=J^*=R$. From \ref{Gor}, $R\simeq Hom_S(I,S)$ as $S$-modules. Now, if $f,g\in \mathfrak{m}$, then it is clear from the free resolution for $I$ over $S$ in \ref{Ifreeres}, that (1) holds.

\par  If $f$ were a unit say, then so is $h_1$, so the resolution in \ref{Ifreeres} is not minimal. In this case, $R$ is Cohen-Macaulay. Therefore, for the proof of (2) and (3) we may assume that $f,g\in \mathfrak{m}$. Since $S/2S$ is regular local (a UFD), we have $Q=(2,zc,ze)$ for some $z\notin 2S$. Then $S/Q$ admits the following free $S$-resolution:
\\
\begin{tikzcd}
\quad \quad \quad \quad \quad \quad \quad 0 \arrow{r} & S \arrow[column sep = large]{r}{[-e,c,-2]^T}
&[5ex] S^3 \arrow{r}{\Phi} & S^3 \arrow{r}{\psi} & S\arrow{r} &S/Q\arrow{r} & 0
\end{tikzcd}
\\
where $\psi(e_1)=2$, $\psi(e_2)=zc$, $\psi(e_3)=ze$ and
  \[
\Phi = \begin{bmatrix} 
    zc & ze &0 \\
    -2 & 0 & e \\
    0  & -2  & -c
    \end{bmatrix}
\]
Note that this is indeed a resolution by the Buchsbaum-Eisenbud criterion \cite{buchsbaum1973makes}[Cor 1].
Thus $p.d_S(R)=p.d._S(Syz_S^2(S/Q))\leq 1$ and (2) holds. Finally, $R$ is Cohen-Macaulay if and only if $c$ or $e$ is a unit. The latter is clearly equivalent to $Q$ being a two generated ideal and thus (3) holds.
\QED
\begin{rem}\normalfont \label{pfg}
Let $S$ be an unramified regular local ring of mixed characteristic two. Note that the ideal $(2,h_1,h_2)\subseteq S$ is a two generated ideal or all of $R$ if and only if the same holds for $(2,f,g)\subseteq S$. 
\end{rem}
\begin{rem}\normalfont \label{R4}
In the context of \ref{R2}, if $Q:=(2,h_1,h_2)S$ is a complete intersection ideal of grade three, then the conductor of $R$ to $A$ is the ideal $I$ in \ref{Ifreeres}.
\par To show this, since the only element of $NNL_1(A)$ is $P=(2,\omega-h_1,\mu-h_2)$ and since the conductor of a ring that satisfies $S_2$ is unmixed, it suffices to show that $I$ is $P$-primary. Let $x\cdot y\in I$ such that $y\in A$ and $x\in A\setminus P$. Certainly $y\in P$, so write $y=2\cdot a_1+(\omega-h_1)\cdot a_2+(\mu-h_2)\cdot a_3$ for some $a_i\in A$.  Since $x\cdot y\in I$ $\text{if and only if}$ $x\cdot (a_2(\omega-h_1)+a_3(\mu-h_2))\in I$, it suffices to show $a_2(\omega-h_1)+a_3(\mu-h_2)\in I$. Lifting to $B:=S[W,U]_{(\mathfrak{m},W-h_1,U-h_2)}$ and denoting lifts by $\sim$, we have for some $\tilde{b_i}\in B$,
\begin{equation}
    \tilde{a_2}\cdot \tilde{x}(W-h_1)+\tilde{a_3}\cdot \tilde{x}(U-h_2)+2\cdot \tilde{b_1}+(WU-h_1h_2)\cdot \tilde{b_2}+(W-h_1)(U-h_2)\cdot \tilde{b_3} \in (F(W),G(U))
\end{equation}
Writing $WU-h_1h_2=(W-h_1)(U-h_2)+h_2(W-h_1)+h_1(U-h_2)$, we have that $\tilde{a_2}\cdot \tilde{x}+\tilde{b_2}\cdot h_2\in \tilde{P}$ where $\tilde{P}:=(2,W-h_1,U-h_2)$. Similarly $\tilde{a_3}\cdot \tilde{x}+\tilde{b_2}\cdot h_1\in \tilde{P}$. Hence $h_1\tilde{a_2}\tilde{x}-h_2\tilde{a_3}\tilde{x}\in \tilde{P}$ and therefore $h_1\tilde{a_2}-h_2\tilde{a_3}\in \tilde{P}$. Since $\tilde{P}+(h_1,h_2)B$ is a grade five complete intersection ideal, $\tilde{a_2}\equiv h_2\cdot z\;\text{mod}\;P$ and $\tilde{a_3}\equiv h_1\cdot z\;\text{mod}\;P$ for some $z\in A$. We have $P\subseteq (I:_A (\omega-h_1,\mu-h_2))$, so \[a_2(\omega-h_1)+a_3(\mu-h_2)\equiv [h_2(\omega-h_1)+h_1(\mu-h_2)]\cdot z\;\text{mod}\;I\] Since $h_2(\omega-h_1)+h_1(\mu-h_2)\in I$, we are done. Thus, $I$ is $P$-primary and hence is the conductor of $R$ to $A$.
\end{rem}

\begin{Example}\label{ex}\normalfont
The conditions in \ref{R2} produce a non-empty class of non Cohen-Macaulay integral closures $R$. In fact they are quite abundant. From $\ref{R2}$, there are two classes of examples, the first one being the case where $Q:=(2,h_1,h_2)$ is grade two with $p.d_S(S/Q)=3$ and the other when $Q$ is grade three perfect. For an example of the first kind, set $S:=\mathbb{Z}[X,Y,V]_{(2,X,Y,V)}$ where $X,Y,V$ are indeterminates over $\mathbb{Z}_{(2)}$ and let \[f=V^2X^2-2X^2+4=(VX)^2+2(2-X^2)\] \[g=V^2Y^2-2Y^2+4=(VY)^2+2(2-Y^2)\] and $\omega^2=f,\mu^2=g$. Then $f,g$ are square free elements that form a regular sequence in $S$. It is straightforward to check that $[L(\omega,\mu):L]=4$. The hypersurface rings $S[\omega]$ and $S[\mu]$ are integrally closed, but the hypersurface ring $S[\omega\mu]$ is not. Since $(2,VX,VY)\subseteq S$ is a grade two ideal such that $p.d_S(S/(2,VX,VY))=3$, by \ref{R2}, $p.d._S(R)=1$.
\par For an example of the second kind, let $S=\mathbb{Z}[X,Y]_{(2,X,Y)}$, where $X,Y$ are indeterminates over $\mathbb{Z}_{(2)}$ and take \[f=-X^2+4=X^2+2(2-X^2)\] \[g=-Y^2+4=Y^2+2(2-Y^2)\]
\QED
\end{Example}

\par We now get to our main theorem showing that $R$ always admits a birational small CM module when $S$ is an unramified regular local ring of mixed characteristic two and $f,g\in S^2$. By \cite{DKPS}, if $S$ is complete with perfect residue field, then this would show that $R$ always admits a small CM module, when $f,g\in S$ are square free and form a regular sequence.
\begin{Lemma}\label{helplemma}
With established notation, $P^*_A$ is generated as an $A$-module by $\{1,\eta\}$, where $\eta:=2^{-1}(\omega+h_1)(\mu+h_2)\in K$.
\end{Lemma}
\demo 
Set $P_1:=(2,\omega-h_1)$ and $P_2:=(2,\mu-h_2)$, so that $P^*=P_1^*\cap P_2^*$. Let $\tilde{P_1}:=(2,W-h_1)\subseteq S[\mu][W]$. It is the maximal minors of 
\[   M = \begin{bmatrix} 
    W-h_1 \\
    2\\

    \end{bmatrix}
    \]
    We have $F(W)\in \tilde{P_1}$, $F(W)=a\cdot (-2)+(W+h_1)(W-h_1)$. Adjoining the appropriate column of coefficients we get
     \[   M' = \begin{bmatrix} 
    W-h_1  & a\\
    2& W+h_1\\
\end{bmatrix}
    \]
From \cite{KU}[Lemma 2.5], $P_1^*$ is generated as $A$-module by $\{M'_{11}/\delta_1,M'_{22}/\delta_2\}$ where $M'_{ii}$ and $\delta_i$ denote the image in $A$ of the $(i,i)$-th cofactor of $M'$ and the $i$-th (signed) minor of $M$ respectively. Therefore, $P_1^*$ is generated as a $A$-module by $\{1,2^{-1}(\omega+h_1)\}$. Identically, $P_2^*$ is generated over $A$ by $\{1,2^{-1}(\mu+h_2)\}$.
Now consider $y\in P^*=P_1^*\cap P_2^*$. Lifting to $B:=S[W,U]$ and denoting lifts by $\sim$
\[2\tilde{y}\in (2,W+h_1)\cap(2,U+h_2)+(F(W),G(U))=(2,(W+h_1)(U+h_2),F(W),G(U))\]
Thus $2y\in (2,(\omega+h_1)(\mu+h_2))A$ and hence this shows $P^*\subseteq A+A\cdot\eta$. The reverse inclusion is clear since $\eta\cdot P\subseteq A$. Thus the proof is complete.
\QED

\begin{Theorem}\label{P11}
Let $S$ be an unramified regular local ring of mixed characteristic two and $f,g\in S^2$.
\begin{enumerate}
\item $R$ is Cohen-Macaulay if and only if one of the following happens
\begin{enumerate}
    \item At least one of $S[\omega],S[\mu]$ is not integrally closed.
    \item $S[\omega],S[\mu]$ are both integrally closed and $fg\notin S^{2\wedge 4}$.
    \item $S[\omega],S[\mu]$ are both integrally closed, $fg\in S^{2\wedge4}$ and $\mathscr{I}:=(2,f,g)\subset S$ is a two generated ideal or all of $R$.
\end{enumerate}
\item If $R$ is not Cohen-Macaulay, $R$ admits a birational small CM module.
\end{enumerate}
\end{Theorem}

\demo
For (1), we have already shown that the conditions in (a) and (b) imply $R$ is Cohen-Macaulay in \ref{T1} and \ref{T1'} respectively. From \ref{R2} and \ref{pfg}, we see that the condition in 1(c) implies that $R$ is Cohen-Macaulay. For the forward implication of (1), the contrapositive follows from \ref{pfg} and $\ref{R2}$. 
\par For (2), by (1) it only remains to be shown that $R$ admits a birational small CM module when $S[\omega], S[\mu]$ are both integrally closed, $fg\in S^{2\wedge 4}$ and $p.d._S(S/Q)>2$. Therefore, assume all of these conditions for the remainder of the proof.

We have from \ref{P10} that for $I:=(2,\omega\mu-h_1h_2)A$, $I^*=R$. Set $M=(IP)^*$, where $P=(2,\omega-h_1,\mu-h_2)$ is the unique height one prime containing $2$ in $A$. Then $(IP)^*$ is an $R$-module since \[(IP)^*=A:_K IP=((A:_K I):P)=(R:_K P)\] We now show that $depth_S(M)=d$. By definition, \[(IP)_A^*=(2\cdot P+(\omega\mu-h_1h_2)\cdot P)_A^*=F_1\cap F_2\] where $F_1=2^{-1}P^*$ and $F_2=(\omega\mu-h_1h_2)^{-1}P^*$. This is because for ideals $J,J'\subseteq A$, $(A:_K J+J')=(A:_K J)\cap (A:_K J')$ as $A$-modules.
\par Now $P$ is $S$-free since $A/P\simeq S/2S$ as $S$-modules and by the depth lemma, $depth_S(P)=d$. By $\ref{Gor}$, $P^*_A\simeq P^*_S$ as $S$-modules, so $P^*$ is Cohen-Macaulay as well and hence $F_1$ and $F_2$ are Cohen-Macaulay. We have the natural short exact sequence of $S$-modules 
\begin{equation*}
\xymatrix@C+1pc{
 0\ar[r] & F_1\cap F_2 \ar[r]
  &
  F_1\oplus F_2  \ar[r]
  &
  F_1+F_2 \ar[r] & 0 
}
\end{equation*}
By the depth lemma, it suffices to show $depth_S(F_1+F_2)\geq d-1$. Clearly $F_1+F_2\simeq F_1'+F_2'$ as $A$-modules and hence $S$-modules where $F_1'=(\omega\mu-h_1h_2)P^*\subseteq A$ and $F_2'=2P^*\subseteq A$. We claim that $F_1'+F_2'=F_2'+(\omega\mu-h_1h_2)$ as ideals of $A$. By \ref{helplemma}, we only need to show that \[v:=2^{-1}(\omega+h_1)(\mu+h_2)(\omega\mu-h_1h_2)\in H:=F_2'+(\omega\mu-h_1h_2)\] Writing \[(\omega\mu-h_1h_2)=(\omega-h_1)(\mu-h_2)+h_2(\omega-h_1)+h_1(\mu-h_2),\] we have 
\begin{equation}\label{4.7E1}
\begin{aligned}
    v&\equiv 2^{-1}(\omega+h_1)(\mu+h_2)(h_2(\omega-h_1)+h_1(\mu-h_2))\; \text{mod}(H)
    \\
    &\equiv ah_2(\mu+h_2)+bh_1(\omega+h_1)\; \text{mod} (H)
    \end{aligned}
    \end{equation} 
   
    Since $S/2S$ is regular local, $h_1\equiv (zc)\text{mod}(2)$, $h_2\equiv (ze)\text{mod}(2)$ for some $z\notin 2S$ and $c,e$ such that $(2,c,e)\subseteq S$ form a regular sequence. From (1) in \ref{T1'}, $fg\in S^{2\wedge 4}$ implies $ah_2^2+bh_1^2\in 2S$ and hence $a-qc^2\in 2S$ and $b+qe^2\in 2S$ for some $q\in S$. Therefore, (\ref{4.7E1}) implies
    \begin{equation}\label{4.7E3}
    \begin{aligned}
    v&\equiv qc^2h_2(\mu+h_2)-qe^2h_1(\omega+h_1)\; \text{mod}(H)
    \\
    &\equiv qce(h_1(\mu+h_2)-h_2(\omega+h_1))\; \text{mod}(H)
    \end{aligned}
    \end{equation} 
    But $\omega\mu-h_1h_2-(\omega-h_1)(\mu-h_2)=h_2(\omega-h_1)+h_1(\mu-h_2)\in H$. Since $2\in H$, $h_1(\mu+h_2)-h_2(\omega+h_1)\in H$ and from (\ref{4.7E3}), $v\in H$. Therefore $F_1'+F_2'=H=(2,(\omega+h_1)(\mu+h_2),\omega\mu-h_1h_2)$ by \ref{helplemma}. From \ref{Ifreeres}, $p.d_S(H)\leq 1$ so that $depth_S(H)\geq d-1$. This completes the proof, and hence $M=(IP)^*$ is a small CM module over $R$. 
\QED
\section*{Acknowledgement}
I would like to thank my Ph.D. advisor Prof. Daniel Katz for suggesting this problem and for his support through the course of this work. I would also like to thank the referee for their careful reading of the manuscript and suggestions for improvement.

\bibliographystyle{alpha}

\end{document}